\documentstyle[a4]{article}

\newcommand{\area}{\mbox{Area}}
\newcommand{\dr}{\partial}

\newcommand{\R}{{\bf R}}
\newcommand{\N}{{\bf N}}
\newcommand{\Z}{{\bf Z}}
\newcommand{\II}{I\hspace{-0.1cm}I}
\newcommand{\III}{I\hspace{-0.1cm}I\hspace{-0.1cm}I}
\newcommand{\tr}{\mbox{tr}}
\newcommand{\rk}{\mbox{rk}}
\newcommand{\ric}{\mbox{ric}}
\newcommand{\deltab}{\overline{\delta}}
\newcommand{\Db}{\overline{D}}
\newcommand{\Kb}{\overline{K}}

\newcommand{\Sb}{\overline{S}}

\newcommand{\ricb}{\overline{\ric}}

\newcommand{\kn}{\bigcirc\hspace{-0.37cm}\wedge \;}

\newcommand{\vs}{\vspace{0.3cm}}

\newtheorem{prop}{Proposition}
\newtheorem{lemma}{Lemma}
\newtheorem{thm}{Theorem}
\newtheorem{cor}{Corollary}
\newtheorem{remark}{Remark}
\newtheorem{ex}{Example}

\makeatletter
\newif\iflabel\labelfalse \let\w@label=\label
\def\label{\global\labeltrue\w@label}
\let\w@eqn=\equation
\def\equation{\global\labelfalse\w@eqn}
\def\endequation{\iflabel\@eeq\else\@eeqw\fi$$\global\@ignoretrue}
\def\@eeq{\eqno \@eqnnum}
\def\@eeqw{\addtocounter{equation}{-1}}
\newskip\bw@\newskip\bws@
\def\dc@@pq{\global \bw@=\belowdisplayskip \global\bws@=\belowdisplayshortskip
   \global\belowdisplayshortskip=3pt plus -3pt
   \global\belowdisplayskip=\belowdisplayshortskip
   \hskip 500pt minus 500pt\relax}
\def\dc@pq{\dc@@pq$$}
\def\fc@pq{$$\hskip 500pt minus 500pt\global
   \belowdisplayshortskip=\bws@\global\belowdisplayskip=\bw@}
\def\coupeq{\dc@pq\fc@pq}
\def\coupepas{\par\noindent\begin{minipage}{\textwidth}}
\def\jusquela{\end{minipage}}
\makeatother

\begin{document}

\title{On the Schl{\"a}fli differential formula}

\author{
Igor Rivin\thanks{Mathematics Department, University of
Manchester, Manchester, England and Mathematics Department, Temple
University; \texttt{rivin@math.temple.edu}} and Jean-Marc
Schlenker\thanks{ Math{\'e}matiques (UMR 8628 du CNRS), B{\^a}t.
425, Uni\-ver\-sit{\'e} de Paris-Sud, F-91405 Orsay Cedex, France.
Currently: FIM, ETHZ, R{\"a}mistrasse 101, CH-8092 Z{\"u}rich,
Switzerland. \texttt{jean-marc.schlenker@math.u-psud.fr}} }

\maketitle

\begin{flushright}
{\it To F. O. Almgren, in memoriam}
\end{flushright}

\begin{quote}

\begin{center} {\bf Abstract} \end{center}

The celebrated formula of Schl{\"a}fli relates the variation of the
dihedral angles of a smooth family of polyhedra in a space form and the
variation of volume. We give a smooth analogue of this classical
formula -- our result relates the variation of the volume bounded by a
hypersurface moving in a general Einstein
manifold and the integral of the variation of the mean curvature.
The argument is direct, and the classical polyhedral result (as well
as results for Lorenzian space forms) is an easy corollary.
We extend it to variations of the metric in a Riemannian Einstein
manifold with boundary.

We apply our results to extend the classical Euclidean inequalities of
Aleksandrov to other $3$-dimensional constant curvature spaces. We
also obtain rigidity results for Ricci-flat
manifolds with umbilic boundaries and existence results  for foliations of
Einstein manifolds by hypersurfaces.

\vspace{0.6cm}

\begin{center} {\bf R{\'e}sum{\'e}} \end{center}

La formule classique de Schl{{\"a}fli} relie la variation des angles
di{\`e}dres d'une famille lisse de poly{\`e}dres dans un espace {\`a} courbure
constante et la variation du volume.

On donne un analogue r{\'e}gulier de cette formule classique -- notre
resultat relie la variation du volume born{\'e} par une hypersurface
se d{\'e}pla{\c c}ant dans une vari{\'e}t{\'e} d'Einstein {\`a}
l'int{\'e}grale de la variation de la courbure moyenne. Puis nous
l'{\'e}tendons aux variations de la m{\'e}trique {\`a} l'int{\'e}rieur
d'une vari{\'e}t{\'e} d'Einstein
riemannienne.

Comme application, on donne un
r{\'e}sultat de rigidit{\'e} pour les vari{\'e}t{\'e}s Ricci-plates {\`a} bord ombilique,
quelques formules concernant les feuilletages de vari{\'e}t{\'e}s d'Einstein par
des hypersurfaces, ainsi que des analogues de la formule
d'Alexandrov-Fenchel dans les vari{\'e}t{\'e}s de dimension 3 {\`a} courbure
constante.

\end{quote}

\vs

Let $M$ be a Riemannian $(m+1)$-dimensional space-form with constant
curvature $K$, and
$(P_t)_{t\in [0,1]}$ a one-parameter family of polyhedra in $M$ bounding
compact domains, all having the same combinatorics. Call $V_t$ the
volume bounded by $P_t$, $\theta_{i,t}$ and $W_{i,t}$ the dihedral angle
and the $(m-1)$-volume of the codimension 2 face $i$ of $P_t$.
The classical Schl{\"a}fli formula (see for instance \cite{milnor-schlafli}
or \cite{geo2}) relates the
variation of $V_t$ and of the angles $\theta_{i,t}$ as follows:

\vs

\newpage
{\bf Schl{\"a}fli's formula:}

\begin{equation} \label{for-schlafli}
\sum_i W_{i,t} \frac{d \theta_{i,t}}{dt} = m K \frac{dV_t}{dt}
\end{equation}
This formula has been extended and used on several occasions recently, see
for instance
\cite{cooriv}, \cite{suarez}, \cite{bonahon}.
The first goal of this
paper is to extend (\ref{for-schlafli}) to deformations of smooth
hypersurfaces in $H^n$. A remarkable point is that the formula which
appears remains valid for deformations of hypersurfaces in Einstein
manifolds. In section 2, we extend it also to a much more general
situation, namely the
deformations of the Einstein metric in an Einstein manifold with
boundary.

When the $W_{i,t}$ are constant, the left-hand side of
(\ref{for-schlafli}) is a
polyhedral analogue of  the variation of the mean curvature integral
of a hypersurface. Indeed, the polyhedral analogue of the mean
curvature integral is $H=\sum W_i \theta_i$ (where $\theta_i$ is the
{\em exterior} dihedral angle at the $i$-th codimension-$2$
face). Using the product rule, we see that
\begin{equation}
\frac{d H}{d t} = \sum \frac{d W_i}{d t} \theta_i + \sum W_i \frac{d
\theta_i}{d t}.
\end{equation}
When the deformation is isometric, the first sum on the right hand
side vanishes, and, combining that with the formula \ref{for-schlafli},
we see that
\begin{equation}
\frac{d H}{d t} = - m K \frac{dV_t}{dt},
\end{equation}
where the minus sign comes from replacing the dihedral angles with
exterior dihedral angles.

When $K=0$, the right-hand side is $0$. This shows that the ``mean
curvature'' of a 1-parameter family of Euclidean polyhedra with
constant induced metric is constant.  This has been used in
\cite{almgren-rivin} to prove, using geometric measure theory methods,
a similar result for isometric deformations of smooth hypersurfaces in
space-forms; we prove this result directly by differential-geometric
methods here (in the more general setting of Einstein manifolds),
indeed, we show the rather strong Theorem \ref{nbhd}. Theorem
\ref{eucl} should be compared with the celebrated mean curvature
variation (under bending) formula of Herglotz; see
\cite{berger-gostiaux}.

Isometric deformations of (smooth) hypersurfaces remain rather
mysterious. It is known since Liebmann \cite{liebmann} (also see
\cite{herglotz}, \cite[vol V]{spivak}) that
(strictly) convex surfaces
in $\R^3$ admit no such deformation, and this has been extended to $S^3$
and $H^3$ by Pogorelov \cite{Po}. On the other hand, it is unknown
whether any smooth closed surface in $\R^3, S^3$ or $H^3$ has a
1-parameter family of isometric deformations, or even if any closed
surface with no open flat region has a smooth infinitesimal isometric
deformation. But it
{\em is} known that there exist non-rigid polyhedra, see
\cite{connelly} and \cite{bleecker}.

A smooth version of
(\ref{for-schlafli}) leads to some rigidity results. Aside from the
remarks following Theorem \ref{eucl}, we prove in Section
3 a rigidity result for Ricci-flat manifolds with umbilic boundaries,
with respect to Einstein deformations of the metric which leave the
induced metric on the boundary fixed.

Section 4 contains some applications of our ``Schl{\"a}fli formula'' to
codimension one foliations of Einstein manifolds. Section 5 recalls
some concepts of integral geometry, to see how some of our results may be
recast in integral-geometric terms. In This is particularly relevant to
section 6 (after giving brief preliminaries in Section 5), where we
extend a classical inequality of Alexandrov for
convex bodies in Euclidean space to the hyperbolic and spherical settings.

This paper can be considered (among other things) as a hint that some
important elements of hyperbolic geometry in dimension 3 can be
extended in higher dimension in the setting of Einstein manifolds
(with negative curvature). Another such hint is given in \cite{ecb},
which contains a partial extension of some classical results of the
theory of convex surfaces in hyperbolic 3-space to Einstein manifolds
with boundary.

The smooth Schl{\"a}fli formula given in section 1 also has
analogous ``higher'' smooth Schl{\"a}fli formulas (see \cite{hsf})
but in constant curvature manifolds only. Going to the polyhedral
case leads to ``higher'' polyhedral Schl{\"a}fli formulas,
relating the variations of the volumes of the $p$-faces to the
variations of the ``curvatures'' of the $(p+2)$-faces. Similar
(polyhedral) formulas were given (in some special cases) in
\cite{suarez}.

Throughout this paper, $M$ is an Einstein manifold of dimension
$m+1\geq 3$, and $D$ is its Levi-Civita connection.
When dealing with a hypersurface $\Sigma$ (resp. with the boundary $\dr
M$), we call $I$ the {\em induced metric}, also called the first
fundamental form, of the corresponding immersion in $M$. $\Db$ is
the Levi-Civita connection of $I$, and $B$ the shape operator defined,
for any $x\in \Sigma$
(resp. $x\in \dr M$) and $X\in T_x\Sigma$ (resp.  $X\in T_x\dr M$)
by
$$ BX = - D_X n~, $$
where $n$ is the oriented normal unit vector to $\Sigma$ (resp. the unit
exterior normal to $\dr M$).
The second fundamental form $\II$ of $\Sigma$ (resp. $\dr M$)
is defined by
$$ \II(X,Y) = I(X, BY)~, $$
and the third fundamental form by
$$ \III(X,Y) = I(BX, BY)~. $$

The trace $H$ of $B$ is called the ``mean curvature'' (some definitions
differ by a factor $m$) and the ``higher mean curvatures''
$H_k, k\geq 1$, are the higher symmetric functions of the principal
curvatures of $\Sigma$ (resp. $\dr M$). For instance, $H_2 = (H^2 -
\tr(B^2))/2$.
$dV, dA$ are the volume elements in $M$ and on $\Sigma$ (resp. $\dr
M$) respectively.

We denote the divergence acting on symmetric tensors by $\delta$, and
its formal adjoint by
$\delta^*$. Therefore, if $h$ is a
symmetric 2-tensor and $(e_i)_{i\in {\N}_{m+1}}$ an orthonormal moving
frame on $M$, then, for any vector $x\in TM$:
$$ (\delta h)(x) = - \sum_i (D_{e_i} h)(e_i, x)~, $$
and, if $v$ is a vector field, $p\in M, x,y\in T_pM$, then:
$$ (\delta^* v)(x,y) = \frac{1}{2} (\langle D_x v, y\rangle + \langle
D_y v, x\rangle )~. $$

We will often implicitly identify (through the metric) vector fields
and 1-forms, as well as quadratic forms and linear morphisms.

\section{Deformation of hypersurfaces}

This section contains an analogue of the Schl{\"a}fli formula for
deformations of (smooth) hypersurfaces in a {\em fixed} Einstein
manifold $M$, which can
be Riemannian or Lorentzian (the other pseudo-Riemannian cases can be
treated in the same way; we have not included them to keep things as
simple as possible). This contrasts with the results in the next
section, where the same formula is proved for variations of the metric
inside a manifold with boundary (which is much more general) but only
when $M$ is Riemannian.

We also show how this ``smooth'' Schl{\"a}fli formula can be used
to recover the classical polyhedral formula (\ref{for-schlafli}),
in the Riemannian and Lorentzian cases.

The techniques here are quite elementary, and use the  method of
moving frames.

\begin{thm} \label{def-hyp-h}
Let $\Sigma$ be a smooth oriented hypersurface in a
(Riemannian) Einstein $(m+1)$-manifold $M$ with
scalar curvature $S$, and $v$ a section of the restriction of $TM$ to
$\Sigma$. $v$ defines a deformation of $\Sigma$ in $M$, which induces
variations $V', H'$ and $I'$ of the volume bounded by $\Sigma$, mean
curvature, and induced metric on $\Sigma$.
Then:
\begin{equation} \label{schlafli-h}
\frac{S}{m+1} V' = \int_{\Sigma} H' + \frac{1}{2} \langle
I',\II\rangle dA
\end{equation}
\end{thm}

$\Sigma$ actually does not need bound a finite volume domain for this
formula to hold. If it doesn't, then $V$ doesn't exist, but its
variation still makes sense (since $\Sigma$ is
homologous to its images under deformations).

\vs

{\bf Proof:} We first prove the formula for $v$ tangent to $\Sigma$, then
we'll check for  normal vector fields. When $v$ is tangent to $\Sigma$,
$V'=0$, and
$$ I'(X,Y)=\langle D_X v,Y\rangle + \langle D_Yv,X\rangle = 2(\delta^*
v)(X,Y) $$
so that:
$$ \int_{\Sigma} \langle I', \II\rangle dA = 2\int_{\Sigma} \langle
\delta^* v, \II\rangle dA = 2 \int_{\Sigma} \langle v,  \delta \II \rangle
dA $$
Let $(e_i)$ be an orthonormal frame for $I$ for which $B$ is diagonal. The
Codazzi equation shows that
$$ (\Db_X \II)(Y,Z) = (\Db_Y \II)(X,Z) + \langle R(X,Y)n, Z\rangle $$
so
\begin{eqnarray}
\langle \delta \II, v\rangle & = & - (\Db_{e_i}\II)(e_i, v)
\nonumber \\
& = & - (\Db_v\II)(e_i, e_i) - \langle R(e_i, v) n, e_i\rangle \nonumber \\
& = & - dH(v) + \mbox{ric}(v,n) \nonumber
\end{eqnarray}
Now $M$ is Einstein and $n$ is orthogonal to $v$, so that:
\begin{equation} \label{del}
\langle \delta \II, v\rangle = -dH(v)
\end{equation}

Therefore:
$$
\int_{\Sigma} \langle I', \II\rangle dA = - 2 \int_{\Sigma} dH(v) dA
$$
This proves the formula when $v$ is tangent to $\Sigma$.

Suppose now that $v$ is a normal vector field, i.e. $v=f n$ for some
function $f$ on $\Sigma$. Since $f$ is the difference between two
strictly positive functions, it is enough to prove the result when $f$
does not vanish. Let $x,y$ be vector fields on $\Sigma$. Choose an
extension of $f n$ to some vector field on a neighborhood $\Omega$ of
$\Sigma$ in
$M$, with $n$ the unit
orthogonal to the image of $\Sigma$ by the flow of $f n$, and $df(n)=0$.
Extend $x, y$
to $\Omega$ by the flow of $f n$, then $[fn, x]=[fn, y]=0$. We now have:
$$
I'(x,y) = fn.\langle x, y \rangle
= \langle D_{fn} x, y\rangle + \langle x, D_{fn} y\rangle
= \langle D_{x} (fn), y\rangle + \langle x, D_{y} (fn)\rangle =
- - -2f \II(x,y)
$$
so $I'=-2 f \II$. One also checks that $D_{fn} n = -Df$, so that:
\begin{eqnarray}
\II'(x,y) & = & - fn.\langle D_x n, y \rangle  \nonumber \\
& = & - \langle D_{fn} D_x n, y\rangle - \langle D_x n, D_{fn} y\rangle
 \nonumber \\
& = & - \langle D_{x} D_{fn} n + R_{fn, x} n + D_{[fn, x]} n, y\rangle -
\langle D_x n, D_{y}(fn) \rangle \nonumber  \\
& = & \langle D_x Df, y\rangle - \langle R_{fn, x} n, y\rangle -
f\III(x,y)  \nonumber
\end{eqnarray}
and
\begin{equation} \label{dII-norm}
\II'= H_f - f \langle R_{n, \cdot}n, \cdot\rangle - f\III
\end{equation}
where $H_f$ is the Hessian of $f$ on $\Sigma$.

Taking the trace of this equation:
$$
H' = \tr(\II') - \langle I', \II\rangle = - \Delta f + f \ric(n,n) - f
\tr(\III) - \langle I', \II\rangle
$$
Now the integral over $\Sigma$ of $\Delta f$ is zero, and the integral of
$f$ is $V'$ because the deformation is normal. The result follows,
because $I'=-2f\II$, so that $-2f\tr(\III)=\langle I', \II\rangle$.
$\Box$

\vs

This formula leads easily to the ``classical'' Schl{\"a}fli formula for
polyhedra in space-forms:

\begin{thm} \label{cor-pol-h}
Let $P$ be a convex polyhedron in a $(m+1)$-dimensional space-form $M$ with
scalar curvature $S$; for any deformation of $P$, the variation $V'$ of
the volume bounded by
$P$ is given in term of the variations $\theta'_i$ of the dihedral angles
at the codimension 2 faces by:
$$ \frac{S}{m+1} V' = \sum_i W_i \theta'_i $$
where $W_i$ is the $(m-1)$-volume of the codimension 2 face $i$.
\end{thm}

{\bf Proof:}
First note that Theorem \ref{def-hyp-h} also applies for deformations of a
$C^{1,1}$, piecewise smooth hypersurface (if the deformation preserves the
decomposition into smooth parts). This is proved by an easy
approximation argument. The formula remains the same, and each term make
sense in this case.

Call $P_\epsilon$ the set of points at distance $\epsilon$ of $P$ on the
outside (i.e. on the side of $P$ which is concave). For $\epsilon$ small
enough, $P_\epsilon$ is a $C^{1,1}$, piecewise smooth hypersurface, and we
can apply
Theorem \ref{def-hyp-h}. Note $I'_\epsilon, \II_\epsilon, H'_\epsilon,
V'_\epsilon$
the quantities corresponding to $I', \II, H', V'$ for $P_\epsilon$. Then:
$$
\frac{S}{m+1} V'_\epsilon =
\int_{P_\epsilon}H'+ \frac{1}{2}\langle I'_\epsilon,\II_\epsilon\rangle
dA~. $$

For $\epsilon$ small enough, we can decompose $P_\epsilon$ as
$$ P_\epsilon = \cup_{k=1}^{m+1} P_{\epsilon, k}~, $$
where $P_{\epsilon, k}$ is the set of points where the normal meets $P$
on a codimension $k$ face. Using the flow of the unit normal vectors to
the $P_\epsilon$, we can also identify $P_\epsilon$ and $P_{\epsilon'}$
for $\epsilon'\neq \epsilon$, so that we can consider e.g. $I'_\epsilon$
as a 1-parameter family of symmetric 2-tensors on a fixed manifold.

If $x\in P_{\epsilon, 2}$, then the normal to $P_\epsilon$ at $x$ meets
some codimension 2 face $F_i$ of $P$; let $\alpha_{i,t}$ be the dihedral
angle at $F_i$. If $v, w\in T_x P_\epsilon$
correspond to vectors orthogonal to $TF$, then
$$ I'_\epsilon(v,w) \simeq \frac{2}{\alpha_{i,t}} \frac{d\alpha_{i,t}}{dt}
I_{\epsilon}(v,w) $$
as $\epsilon\rightarrow 0$. On the other hand,
$$ \II_\epsilon(v,w) \simeq \frac{1}{\epsilon} I_{\epsilon}(v,w) $$
If $v,w$ now correspond to vectors in $TF_i$, then
$$ I'_\epsilon(v,w) = O(1)~, $$
while
$$ \II_{\epsilon}(v, w) = O(\epsilon)~. $$
Using those 2 cases, we see that, at any point in $P_{\epsilon, 2}$:
$$ \langle I'_\epsilon, \II_\epsilon\rangle \simeq \frac{2}{\epsilon
\alpha_{i,t}} \frac{d\alpha_{i,t}}{dt}~. $$
Now the volume element of $P_{\epsilon, 2}$ is equivalent to $\epsilon$ as
$\epsilon\rightarrow 0$, so:
$$ \lim_{\epsilon\rightarrow 0} \int_{P_\epsilon,2} \langle I'_\epsilon,
\II_\epsilon\rangle dA = \sum_i 2 W_{i,t} \frac{d\alpha_{i,t}}{dt}~. $$

For $P_{\epsilon, 1}$ (that is, for codimension 1 faces), only vectors
parallel to the faces have to be taken into account, and their
contribution is of order $O(\epsilon)$ (as above for $P_{\epsilon, 2}$).

For $P_{\epsilon, k}$ with $k\geq 3$, the same reasoning shows that only
vectors orthogonal to the faces count; if $v, w$ are such vectors, then
$$ I'_{\epsilon}(v,w) = O(I_{\epsilon}(v,w))~, $$
while
$$ \II_{\epsilon}(v,w) = O(I_{\epsilon}(v,w)/\epsilon)~, $$
and the volume element on $P_{\epsilon, k}$ is as $O(\epsilon^{k-1})$,
so
$$ \lim_{\epsilon\rightarrow 0} \int_{P_{\epsilon, k}} \langle
I'_{\epsilon}, \II_{\epsilon}\rangle dA = 0~. $$
It is also easy to check that
$$
\lim_{\epsilon\rightarrow 0} \int_{P_\epsilon} H'_\epsilon dA = 0~,
$$
and this leads to the formula.
$\Box$

\vs

Of course, $P$ does not need to be convex: once the corollary is proved
for convex polyhedra, it is clear that it also applies to
non-convex ones, since they can be decomposed into convex pieces.

The proof of Theorem \ref{def-hyp-h} also applies to the Lorentzian case. The
only difference is
that now $g(n,n)=-1$, so the volume variation has a minus sign in the
formula.

\begin{thm}
Let $\Sigma$ be a smooth oriented space-like hypersurface in a
Lorentzian Einstein $(m+1)$-manifold $(M, g)$ with
$\ric_g = mkg$, and let $v$ be a section of the restriction of $TM$ to
$\Sigma$. $v$ defines a deformation of $\Sigma$ in $M$, which induces
variations $V', H'$ and $I'$ of the volume bounded by $\Sigma$, mean
curvature, and induced metric on $\Sigma$.
Then:
\begin{equation} \label{schlafli-ds}
- - - m k V' = \int_{\Sigma} H'+ \frac{1}{2}\langle
I',\II\rangle dA
\end{equation}
\end{thm}

Here again, the
volume might be defined only up to an additive constant (for instance as
the volume bounded by $\Sigma$ and some fixed homologous hypersurface
$\Sigma_0$), but its
variation is well defined. For instance, if $\Sigma$ is a compact space-like
hypersurface in the de Sitter space, its ``volume'' can be defined as
the oriented volume of the domain bounded by $\Sigma$ and by some
space-like totally geodesic hyperplane $S_0$. This volume actually does not
depend on $S_0$, because if $S_1$ is some other totally geodesic
hyperplane, then, as (\ref{schlafli-ds}) shows, the oriented volume of
the domain bounded by $S_0$ and $S_1$ is zero.

This lemma could actually be extended almost without change to other
pseudo-Riemannian manifolds, and also to hypersurfaces which are not
space-like.

Applying this lemma to the set of points at distance $\epsilon$ from a
polyhedron in $S^n_1$ (as above in Theorem \ref{cor-pol-h}), one
obtains the Schl{\"a}fli formula for de Sitter polyhedra as in
\cite{suarez} (where it was proved for simplices using a more
combinatorial approach).

\begin{thm} \label{cor-pol-dS}
Let $P$ be a convex space-like polyhedron in the de Sitter space
$S^{m+1}_1$, which is dual to a hyperbolic polyhedron.
For any deformation of $P$, the variation $V'$ of the volume bounded by
$P$ is given in term of the variations $\theta'_i$ of the dihedral angles
at the codimension 2 faces by:
$$ m V' + \sum_i W_i \theta'_i = 0 $$
where $W_i$ is the $(m-1)$-volume of the codimension 2 face $i$.
\end{thm}

The conditions that $P$ is convex and dual to a hyperbolic polyhedron
are actually
not necessary, and the formula even remains valid for many polyhedra that are
not space-like. It then helps to use a definition of angles and volume
well adapted to this Lorentzian setting, i.e. with complex values (as in
\cite{shu}). It is not obvious how to give a complete proof using
smooth formulas (as above) but many cases can be treated simply by using
sums or differences of polyhedra for which smooth formulas work. For
instance, this can be done for all space-like polyhedra.

Theorem \ref{def-hyp-h} applied in Euclidean space leads to the

\begin{thm} \label{eucl}
In $\R^{m+1}$, the integral of
the mean curvature remains constant under an isometric deformation of a
compact hypersurface.
\end{thm}

On the other hand, the
integral mean curvature is {\bf not} determined by the metric on $\dr
M$: this is already visible in $\R^3$. Namely,
some metrics on $S^2$ admit two isometric embeddings in $\R^3$: the
classical example is that a (topological) sphere in $\R^3$ which is
tangent to a plane along a circle can be ``flipped'' so as to obtain
another embedding with the same induced metric \cite{spivak}. Those two
embeddings do not in general have the same integral mean curvature --
and thus we have a complicated way of seeing that the two flipped
surfaces cannot be bent one into the other.

The analogue of Theorem \ref{eucl} is also true, but in a pointwise
sense, for the higher mean curvatures:

\begin{thm} \label{higher}
In $\R^{m+1}$, the integral of
$H_k$ ($k\geq 2$) remains constant in an isometric deformation of a
hypersurface.
\end{thm}

This comes from the following (probably classical) description
of the possible
isometric deformations of a hypersurface for $m+1\geq 4$:

\begin{remark} \label{deform}
Let $(\Sigma_t)_{t\in [0,1]}$ be a 1-parameter family of hypersurfaces
in a space-form, such
that the induced metric $I_t$ is constant to the first order at
$t=0$. Then, at each point, one of the following is true:
\begin{itemize}
\item $\II_0=0$;
\item $\rk(\II_0)\leq 2$, and $\II_0'$ vanishes on the kernel
of $\II_0$;
\item $\II_0'=0$.
\end{itemize}
where $\II_t$ is the second fundamental form of $\Sigma_t$, and $\II'_t$
its variation.
\end{remark}

Theorem \ref{higher} clearly follows, because $H_k'$ is zero
for $k\geq 3$ in each case, and the Gauss formula gives the proof for
$k=2$.

\vs

{\bf Proof of Remark \ref{deform}:} Choose an orthonormal frame
$(e_1, \cdots, e_{m})$ on $\Sigma_0$ for which $\II_0$ is
diagonal, with eigenvalues $(k_1, \cdots, k_m)$. By the Gauss
formula, $\II_t \kn \II_t$ (where $\kn$ is the Kulkarni-Nomizu
product) is determined by the induced metric, and is thus
independent on $t$. Therefore, for any choice of indices
$p,q,r,s$: $$ \II_0(e_p, e_s) \II'_0(e_q, e_r) + \II'_0(e_p, e_s)
\II_0(e_q, e_r) = \II_0(e_p, e_r) \II'_0(e_q, e_s) + \II'_0(e_p,
e_r) \II_0(e_q, e_s) $$ Taking $p,q,r$ distinct but $s=p$ shows
that
\begin{equation} \label{eg1}
k_p \II'_0(e_q, e_r) = 0
\end{equation}
while taking $p=s\neq q=r$ leads to:
\begin{equation} \label{eg2}
k_p \II'_0(e_q, e_q) + k_q \II'_0(e_p, e_p) = 0
\end{equation}

Consider the case where $\rk(\II_0)\geq 3$. For each choice of
$p,q,r$ with $k_p, k_q, k_r\neq 0$, adding eq.~(\ref{eg2})
(divided by $k_p k_q$) for the pairs $(p,q)$ and $(p,r)$ and
subtracting the same equation for the pair $(q,r)$ shows that
$\II'_0(e_p, e_p)=0$, and the same for $q,r$, so we already see
that all diagonal terms of $\II'_0$ are zero. Then eq.~(\ref{eg1})
shows that all non-diagonal terms are zero too, so $\II'_0=0$.

If $\rk(\II_0)\leq 2$ but $\II_0\neq 0$, then eq.~(\ref{eg1}) and
eq.~(\ref{eg2}) easily show that $\II'_0=0$ except maybe in the
subspace generated by the eigenvectors of $\II_0$ with non-zero
eigenvalue. $\Box$

Theorems \ref{eucl} and \ref{higher} can be combined to give the
following geometric statement. Denote by $\Sigma^\epsilon_t$ the
parallel surface at distance $\epsilon$ from $\Sigma_t$. It is
well-known (see, eg, Santalo's book \cite{santalo}) that the area
of $\Sigma^\epsilon$ is a polynomial in $\epsilon$ where the
coefficient of $\epsilon^k$ is (essentially) the $k$-th mean
curvature of $\Sigma$. The two Theorems $\ref{eucl}$ and
$\ref{higher}$ can than be combined as stating that:

\begin{thm}\label{nbhd}
The area of $\Sigma_t^\epsilon$ stays constant when $\Sigma_t$ is a
bending of $\Sigma_0.$
\end{thm}

\section{Einstein manifolds with boundary}

In this section, $(M, \dr M)$ is a compact manifold with boundary with
an Einstein metric $g$ of scalar curvature $S$. We will
prove the same formula as in the previous section, but in a much more
general setting: instead of moving a hypersurface in an Einstein
manifold, we will be changing the metric (among Einstein metrics of
given scalar curvature) inside this  manifold with boundary. Although
the two
operations are equivalent in dimension at most 3, moving the inside
metric is much more general in higher dimension.
On the other hand, our proof only works for Riemannian Einstein
manifolds. It is not obvious whether it can be extended to the
pseudo-riemannian setting.

As always when studying deformations of Riemannian metrics, we need
put some kind of restriction to remove the indeterminacy coming from the
fact that some deformations are geometrically trivial, that is, they
just correspond to the action of vector fields on the metric. We
prevent those deformations in the same way as e.g. in \cite{graham-lee},
\cite{deturck}
or \cite{biquard}, that is, we only consider metric variations $h$ such
that $2\delta h + d\tr h = 0$. The following proposition shows that we
don't forget any metric variation when doing this.

\begin{prop}
Let $h'$ be a smooth variation of $g$. Suppose that either $S\leq 0$, or
that $M$ is strictly convex. There exists another smooth
variation $h$ of $g$ such that $2 \delta h + d\tr(h)=0$ and that
$h=h'+\delta^* v_0$, where $v_0$ is a vector field vanishing on $\dr M$.
\end{prop}

{\bf Proof:}
Suppose $v$ is a vector field on $M$, let $h=h'+\delta^* v$. Then
$$ 2\delta h  + d\tr h  = 2\delta h' + d\tr h' + 2 \delta(\delta^*v) +
d\tr(\delta^*v)~. $$
Now, if $x$ is a vector field on $M$:
\begin{eqnarray}
2 \delta(\delta^*v)(x) & = & - \sum_i 2 (D_{e_i}(\delta^* v))(e_i, x)
\nonumber \\
& = & \sum_i - e_i.(2 \delta^*v)(e_i, x) + (2 \delta^*v)(D_{e_i}e_i, x)
+ (2 \delta^*v)(e_i, D_{e_i}x)~, \nonumber
\end{eqnarray}
so
$$ 2 \delta(\delta^*v)(x) = \sum_i - e_i.(\langle D_{e_i} v, x\rangle +
\langle D_x v,
e_i\rangle) + \coupeq + \langle D_{D_{e_i}e_i}v, x\rangle + \langle D_{x} v,
D_{e_i} e_i \rangle +  \langle D_{e_i}v, D_{e_i}x\rangle + \langle
D_{D_{e_i}x}v, e_i\rangle $$
and
$$ 2 \delta(\delta^*v)(x) = \sum_i \langle - D_{e_i} D_{e_i} v +
D_{D_{e_i}e_i}v, x\rangle +
\langle - D_{e_i}D_x v + D_{D_{e_i}x}v, e_i\rangle~.
$$
On the other hand:
$$
d(\tr(\delta^* v))(x)
= x.\left(\sum_i \langle D_{e_i}v, e_i\rangle \right)
= \sum_i \langle D_x D_{e_i}v, e_i\rangle + \langle D_{e_i}v,
D_{x} e_i\rangle~,
$$
so that
$$ (2 \delta(\delta^*v) + d\tr(\delta^*v))(x) = \sum_i \langle
- - - D_{e_i} D_{e_i} v + D_{D_{e_i}e_i} v, x\rangle - \langle
R_{e_i, x}v,
e_i\rangle + \langle D_{D_xe_i}v, e_i\rangle + \langle D_{e_i}v, D_x
e_i\rangle~. $$
If $\omega$ is the connection form of the frame $(e_i)_{i\in
\N_{m+1}}$, then
$$ \sum_i \langle D_{D_xe_i}v, e_i\rangle + \langle D_{e_i}v, D_x
e_i\rangle=\sum_i (2 \delta^*v)(e_i, D_x e_i) = \langle 2 \delta^*v,
\omega(x)\rangle = 0 $$
because $\delta^* v$ is symmetric and $\omega(x)$ is
skew-symmetric. Therefore,
$$ (2 \delta(\delta^*v) +d\tr(\delta^*v))(x) =
\langle D^* D v, x\rangle - \ric(v,x)
= \langle D^* D v, x\rangle - \frac{S}{m+1} \langle v, x\rangle~, $$
so
\begin{equation} \label{wb-vectors}
2 \delta(\delta^*v) +d\tr(\delta^*v) = D^* D v - \frac{S}{m+1} v~.
\end{equation}

To prove the proposition, we have to solve the elliptic problem:
\begin{equation} \label{elliptic}
\left\{
\begin{array}{rcl}
D^* D v - \frac{S}{m+1} v & = & - (2\delta h' + d\tr(h')) \\
v_{|\dr M} & = & 0
\end{array}
\right.
\end{equation}

Call $\Gamma^1_0 TM$ the space of vector fields on $M$ which are in the
Sobolev space $H^1$ and whose trace on $\dr M$ vanishes (this
essentially means that they are zero on $\dr M$), and define
$$
\begin{array}{rrcl}
F: & \Gamma^1_0 TM & \rightarrow & \R \\
& v & \mapsto & \frac{1}{2}\int_M \langle Dv, Dv\rangle - \frac{S}{m+1}
\langle v,v\rangle dV + \int_M \langle 2\delta h' + d \tr(h'),
v\rangle dV~.
\end{array}
$$
Then $F$ is strictly convex, and moreover it is coercive; this is clear if
$S<0$, and, if $S=0$, it follows from the
Poincar{\'e} inequality for vector fields vanishing on $\dr M$:
$$ \exists C, \forall v\in \Gamma^1_0 TM, \int_M \langle Dv, Dv\rangle dV
\geq C \int_M \langle v,v\rangle dV $$
If $S>0$, a more careful argument is necessary. Let $u:=\| v\|$. Then
$$ \langle Du, Du\rangle \leq \langle Dv, Dv\rangle $$
and
$$ \int_M \langle Du, Du\rangle dV \geq \lambda_1 \int_M u^2 dV~, $$
where $\lambda_1$ is the first eigenvalue of the Dirichlet problem for
the Laplacian on $M$. But it is known (see \cite{reilly}, \cite{kasue})
that, for $M$ convex, and under the hypothesis that the Ricci curvature
is bounded below by $S/(m+1)$
$$ \lambda_1 \geq \frac{S}{m}~, $$
with equality if and only if $M$ is a hemisphere. Therefore,
$$ \int_M \langle Dv, Dv\rangle dV \geq \frac{S}{m}\int_M \langle v,
v\rangle dV $$
and $F$ is again coercive.

Therefore, $F$ admits a unique minimum $v_0$ on $\Gamma^1_0 TM$, which
is smooth by standard elliptic arguments. Then, for all $u\in \Gamma TM$,
$$ (T_{v_0}F)(u) = 0 $$
so that
$$ \int_M \langle Dv_0, Du\rangle - \frac{S}{m+1} \langle v_0, u\rangle
dV + \int_M \langle 2\delta h' + d\tr(h'), u \rangle dV = \coupeq =
\int_M \langle D^* D v_0 - \frac{S}{m+1} v_0 + 2\delta h' + d\tr(h'),
u\rangle dV = 0 $$
and $ D^* D v_0 - \frac{S}{m+1} v_0  = - 2\delta h' - d\tr(h')$ as
needed.
$\Box$

\vs

Another way of solving eq.~(\ref{elliptic}) would be to check that
it has index $0$, and that if $ 2 \delta (\delta^* v) +
d\tr(\delta^* v) = 0$ on $M$ and $v=0$ on $\dr M$, then $v\equiv
0$.

If $g$ is an Einstein metric, we say that a 2-tensor $h$ is an
``Einstein variation'' of $g$ if the associated variation of the metric
induces a variation of the Ricci tensor which is proportional to $h$,
so that $g+\epsilon h$ remains, to the first order, an Einstein
manifold with constant scalar curvature.

\begin{thm} \label{schlafli-einstein}
Let $h$ be a smooth Einstein variation of $g$. Then:
\begin{equation} \label{schlafli-e}
\frac{S}{m+1} V' = \int_{\dr M} H' + \frac{1}{2} \langle
h_{|\dr M},\II\rangle dA
\end{equation}
\end{thm}

{\bf Proof:} By the previous proposition, we can suppose that $2\delta
h + d\tr(h)=0$. First, we compute the variation $\II'$ of $\II$ on
$\dr M$. Let
$x$ be a vector field on $M$ so that $D_n x = 0$. Extend
$n$ to a unit vector field on $M$ such that $D_n n=0$. Then
$$ 2 \II(x,x) = - 2 \langle D_x n, x\rangle = - n.\langle x,x \rangle
- - - 2 \langle[x,n], x\rangle~. $$
Now, since $n$ remains the unit normal to $\dr M$
$$ n'= - \frac{\tau}{2}n - a $$
where $\tau = h(n,n)$ and $a\in T\dr M$ is such that for any vector
$y\in T\dr M$, $\langle
y, a\rangle = h(n,y)$. Therefore,
\begin{eqnarray}
2 \II'(x,x) & = & - n.h(x,x) - 2h([x,n],x) +
\left(\frac{\tau}{2}n+a\right).\langle
x,x\rangle + \langle[x, \tau n + 2a], x\rangle \nonumber \\
& = & - (D_n h)(x,x) + 2h(Bx, x) + a.\langle x,x \rangle - \tau \langle
Bx, x\rangle + 2\langle [x,a], x\rangle \nonumber \\
& = & - (D_n h)(x,x) + 2h(Bx, x) +
2\langle D_x a, x\rangle - \tau \langle Bx, x\rangle~. \nonumber
\end{eqnarray}
To go further, we note $\deltab$ the divergence on $\dr M$, $\alpha$ the
$1$-form dual to $a$ on $\dr M$, and $t:=\tr(h)$.

If $(u_1, \cdots, u_n)$ is an orthonormal frame on $\dr M$ for which
$\II$ is diagonal, extended on $M$ so that $D_n u_i=0$, we have:
$$ 2 \tr(\II') = - \sum_i (D_n h)(u_i, u_i) + 2 \langle h, \II\rangle
- - - 2 (\deltab \alpha) - \tau \tr(\II)~. $$
But
\begin{eqnarray}
- - - \sum_i (D_n h)(u_i, u_i) & = & - dt(n) + (D_n h)(n,n) \nonumber \\
& = & - dt(n) - (\delta h)(n) - \sum_i (D_{u_i} h)(u_i, n) \nonumber \\
& = & - \frac{dt(n)}{2} - \sum_i (D_{u_i} h)(u_i, n) \nonumber \\
& = & - \frac{dt(n)}{2} - \sum_i u_i.\alpha(u_i) + h(D_{u_i}
u_i, n) + h(u_i, D_{u_i}n) \nonumber \\
& = & - \frac{dt(n)}{2} - \sum_i (\Db_{u_i}\alpha)(u_i) +
\II(u_i, u_i) \tau - h(u_i, B u_i)  \nonumber
\end{eqnarray}
and, finally,
\begin{equation} \label{II}
2 \tr(\II') = - \frac{dt(n)}{2} + \langle h, \II\rangle - (\deltab
\alpha)~.
\end{equation}
The statement that $h$ is an Einstein variation of $g$ can be written (see
\cite{Be},
chapter 12) as the equation:
$$ D^* Dh - 2\overline{R}h - 2\delta^* \delta
h - Ddt = 0~, $$
where $\overline{R}$ is the curvature operator acting on symmetric
2-tensors; and, since $2\delta h+dt=0$:
$$ D^* Dh - 2 \overline{R}h =0~. $$
Taking the trace of this equation, we find that:
$$ \Delta t = \frac{2 S}{m+1} t~. $$

An elementary computation shows that the variation of the volume
of $M$ is equal to half the integral of the trace $t$ of $h$: $$ 2
V' = \int_M t dV~. $$ But $$ \frac{2 S}{m+1} \int_M t dV = \int_M
\Delta t dV = - \int_{\dr M} dt(n) dA $$ and, using
eq.~(\ref{II}), we obtain: $$ 2\int_{\dr M} H' dA = \int_{\dr M}
2\tr(\II') - 2 \langle h_{|\dr M}, \II\rangle dA = \coupeq = -
\int_{\dr M} \frac{dt(n)}{2} + \langle h_{|\dr M}, \II \rangle +
(\deltab \alpha) dA = \frac{2 S}{m+1} V' - \int_{\dr M} \langle
h_{|\dr M}, \II \rangle dA $$ from which the result follows.
$\Box$

\vs

Formula (\ref{schlafli-e}) is even simpler for variations which
vanish on $\dr M$:

\begin{thm} \label{cor-e}
If $h$ is a smooth Einstein variation of $g$ which does not change the
induced metric on $\dr M$, then
$$ \int_{\dr M} H' dV= \frac{S}{m+1} V' $$
\end{thm}

In particular, for $S=0$, this implies that the integral of the mean
curvature of the boundary is constant under an Einstein variation which
does not change the induced metric on $\dr M$; this is a direct
generalization of Theorem \ref{eucl}.

A more interesting application can be found by looking at
``singular objects'', just as we did to get the polyhedral Theorem
\ref{cor-pol-h} from the smooth Lemma \ref{def-hyp-h}. There are
no polyhedra in general Einstein manifolds, but we can check what
happens when we deform Einstein manifolds with cone singularities.
It should be pointed out that some of the most interesting modern
uses of the classical Schl{\"a}fli formula concern hyperbolic
3-dimensional cone-manifolds.

Let $M$ be a compact $(m+1)$-manifold, and $N$ a compact
codimension 2 submanifold of $M$. Suppose $(g_t)$ is a 1-parameter
family of Einstein metrics with fixed scalar curvature $S\leq 0$ on
$M\setminus N$, with a conical singularity
on $N$ in the sense that, in normal coordinates around $N$, $g_t$ has an
expansion like:
$$ g_t = h_t + dr^2 + r^2 d\theta^2 + o(r^2) $$
where $h_t$ is the metric induced on $N$ by $g_t$, and $\theta\in
\R/\alpha_t \Z$ for some $\alpha_t\in \R$. Call $V_t$ the volume of
$(M\setminus N, g_t)$, and $W_t$ the volume of $(N, h_t)$. Then:

\begin{cor}
$V_t$
varies as follows:
$$ \frac{S}{m+1} \frac{d V_t}{dt} = W_t \frac{d \alpha_t}{dt}~. $$
The same formula of course remains true if $N$ has several connected
components, each with a different value of $\alpha_t$.
\end{cor}

{\bf Proof:} It is similar to the proof of Theorem \ref{cor-pol-h}, so
we go a little faster. Set
$$ N_\epsilon(t) = \{ x\in M \; , \; d(x, N)\geq \epsilon\alpha_t \} $$
Apply Theorem \ref{schlafli-einstein} to the boundary of $N_\epsilon(t)$
and take the limit as
$\epsilon\rightarrow 0$. Again, we consider $\dr N_{\epsilon}$ as a
fixed manifold (diffeomorphic to $N\times S^1$) with a 1-parameter
family of metrics depending on $\epsilon$. If $v, w\in T\dr N_{\epsilon}$
correspond to vectors in $TN$, then
$$ I'_{\epsilon}(v,w) = O\left( 1\right) $$
while
$$ \II_{\epsilon}(v,w) = O(\epsilon) $$
If $v\in T\dr N_{\epsilon}$ corresponds to a vector normal to $N$, then
$$ I'_{\epsilon}(v,v) \simeq \frac{2}{\alpha_t} \frac{d\alpha_t}{dt}
I_{\epsilon}(v,v) $$
while
$$ \II_{\epsilon}(v,v) \simeq \frac{I_{\epsilon}(v,v)}{\epsilon} $$
so that finally
$$ \langle I'_\epsilon, \II_\epsilon\rangle \simeq \frac{2}{\epsilon
\alpha_t} \frac{d\alpha_t}{dt} $$
On the
other hand, the mean curvature of the boundary (given only by the
$\dr/\dr \theta$ direction) is
$$ H_{\epsilon} = \frac{1}{\epsilon} + o\left(\frac{1}{\epsilon}\right) $$
and
$$ H'_{\epsilon}= o\left(\frac{1}{\epsilon}\right)$$
so that we finally find that
$$ \frac{S}{m+1} V_t' = \lim_{\epsilon\rightarrow 0} \int_{\dr
N_\epsilon(t)} \frac{1}{2}\langle I'_{\epsilon}, \II_{\epsilon}\rangle
dA = \frac{d\alpha_t}{dt} W_t $$
$\Box$

\vs

{\bf Note:} The same computation could be made with $N$ replaced by a
stratified submanifold; the same result follows.

{\bf Example:} take $m+1=3$ in the previous example. We find the
Schl{\"a}fli formula for the variation of the volume of a hyperbolic
cone-manifold.

\vs

\section{Applications to rigidity}

In this section, we use the Schl{\"a}fli formula above to prove a rigidity
result for Ricci-flat manifolds with umbilic boundary; it is a
generalization of the classical result (see \cite{spivak}) that the
round sphere is rigid in
$\R^3$, that is, it can not be deformed smoothly without changing its induced
metric.

This kind of rigidity result could be used in the future to prove
that, given a Ricci-flat manifold $M$ with umbilic boundary and
induced metric $g_0$ on the boundary, any metric close to $g_0$ on
$\dr M$ can be realized as induced on $\dr M$ by some Ricci-flat
metric on $M$. In this setting, rigidity corresponds to the local
injectivity of an operator sending the metrics on $M$ to the
metrics on $\dr M$. In dimension 3, this would be a part of the
classical result (see \cite{N}) that metrics with curvature $K>0$
on $S^2$ can be realized as induced by immersions into $\R^3$.
This circle of ideas is illustrated in \cite{ecb}. It is rather
remarkable that the same condition (that the boundary is umbilic)
appears both here and in \cite{ecb}, in the same kind of rigidity
questions, but in a very different way.

The first point is to understand what an umbilic hypersurface in an
Einstein manifold is. By definition, if $N$ is a Riemannian manifold and is
$S$ a hypersurface, then $S$ is umbilic if, at each point $s\in S$, $\II$
is proportional to $I$, with a proportionality constant $\lambda(s)$
depending on $s$. Now

\begin{remark}
If $N$ is Einstein, then $\lambda$ is constant on each connected
component of $S$.
\end{remark}

{\bf Proof: }
Let $B$ be the shape operator of $S$, and $n$ the unit normal. For $s\in
S$ and $x,y\in T_s S$, the Codazzi formula asserts that:
$$ (d^{\Db} B)(x,y) = R_{x,y} n $$
where $\Db$ is the Levi-Civita connection of $S$, and $R$ the curvature
operator of $N$. Since $\II = \lambda I$, this means that:
$$ (d\lambda (x)) y - (d\lambda(y)) x = R_{x,y} n $$
Let $(e_i)_{1\leq i\leq m}$ be a moving frame on $S$. Taking the trace of
the previous expression with respect to $y$ shows that:
$$ m d\lambda(x) - \sum_{i=1}^m d\lambda(e_i) \langle x, e_i\rangle =
- - - \ric(x, n) $$
and $\ric(x,n)=0$ because $N$ is Einstein and $n$ is orthogonal to
$x$. Therefore:
$$ (m-1) d\lambda(x) = 0 $$
for any tangent vector $s$ to $S$, so $\lambda$ is locally constant.
$\Box$

\begin{cor} \label{umbi}
Umbilic hypersurfaces of Einstein manifolds are analytic.
\end{cor}

A rather clumsy way to prove this is to note that, because of the
previous remark, umbilic hypersurfaces are locally graphs of solutions
of some elliptic PDE with analytic coefficients (because Einstein
metrics are analytic, see \cite{Be}). A classical elliptic smoothness
theorem then gives the result. As a consequence:

\begin{cor} \label{extension}
Let $(\Sigma, h)$ be a compact analytic Riemannian $m$-manifold, and
$\lambda, S\in \R$. There exists at most one germ of Einstein $(m+1)$-manifold
around $\Sigma$ with scalar curvature $S$ which induces $h$ and for
which $\Sigma$ is umbilic with $\II = \lambda I$.
\end{cor}

{\bf Proof: } Let $M$ be such a germ of Einstein manifold around
$\Sigma$, and $g$ its metric. By taking the geodesic flow of the
exponential normal to $\Sigma$ in $M$, we see that $g$ can be
locally written, in some neighborhood $V$ of $\Sigma$, as $$ g =
k_t + dt^2~. $$ We call $\II_t$ the second fundamental form of the
hypersurface $\Sigma \times \{ t\}$ for $g$, and $\III_t$ the
corresponding third fundamental form. Choose $m\in \Sigma$, and
$x,y\in T_m \Sigma$. Then a classical computation (which was done,
in a slightly more general case, in the proof of Lemma
\ref{def-hyp-h}) shows that: $$ \frac{d\II_t}{dt}(x,y) = -
\III_t(x,y) + \langle R(x,n)y, n \rangle~. $$ Let $(e_i)$ be an
orthonormal frame at $m$. Call $R$ the Riemann curvature tensor of
$g$, and $R_t$ the curvature tensor of $k_t$. By the Gauss
formula, for $i\in \N$: $$ \langle R_t(x, e_i)y, e_i\rangle =
\langle R(x, e_i)y, e_i\rangle + (\II_t \kn \II_t)(x, e_i, y,
e_i)~, $$ where $ \kn $ is the Kulkarni-Nomizu product (see
\cite{Be}). Taking the trace and calling $\ric$ the Ricci
curvature of $M$ and $\ric_t$ the Ricci curvature of $k_t$ leads
to: $$ \ric_t(x,y) = \ric(x,y) - \langle R(x,n)y, n\rangle +
\sum_i (\II_t \kn \II_t)(x, e_i, y, e_i) $$ and $$ \sum_{i=1}^m
(\II_t \kn \II_t)(x, e_i, y, e_i) = H_t \II_t(x,y) - \III_t(x,y)
$$ so that $$  \frac{d\II_t}{dt} = \ric - \ric_t + H_t \II_t - 2
\III_t~. $$ Now $\II_t = - dk_t/dt$, and $\ric_t$ can be
considered as a second-order elliptic operator in $k_t$. So $k_t$
satisfies
\begin{equation} \label{edp}
\frac{d^2 k_t}{dt^2} = P(k_t) + R\left(\frac{k_t}{dt}\right)
\end{equation}
where $P$ is a second-order elliptic operator and $R$ is an
operator of degree 0 (all solutions of eq.~(\ref{edp}) do not
correspond to germs of Einstein manifolds around an umbilic
surface; for instance, for $S=0$ and $m=2$, only constant
curvature metrics can be induced on umbilic surfaces in Euclidean
space).

Now we can apply the Cauchy-Kowalevskaya theorem (or the
Cartan-K{\"a}hler theorem) which shows that eq.~(\ref{edp}) has a
unique analytic solution. On the other hand, Corollary \ref{umbi}
asserts that $\Sigma$ has to be analytic in $M$, and then $k_t$
has to be analytic, and also $g$ (see \cite{Be}, 5.F). Equation
(\ref{edp}) therefore has a unique solution, which is analytic.
$\Box$

\vs

Note that the solutions of equation (\ref{edp}) might not
correspond to a germ of Einstein manifold around $\Sigma$, but
such a germ is always obtained as a solution of eq.~(\ref{edp}),
and so is unique.

The next step is an inequality concerning the integral of the mean
curvature squared.

\begin{prop} \label{umbilic}
Let $(M, g)$ be an Einstein manifold with boundary, with scalar
curvature $S$. Call $\Sb$ the scalar curvature of $(\dr M, g_{|\dr M})$.
Then the mean curvature $H=\tr(\II)$ of $\dr M$ satisfies
$$ \frac{\Sb}{m-1} - \frac{S}{m+1} \leq \frac{H^2}{m}~, $$
with equality if and only if $\dr M$ is umbilic.
\end{prop}

{\bf Proof:} Call $\ricb$ the Ricci curvature of $(\dr M, g_{|\dr M})$,
and $K(x,y)$ (resp. $\Kb(x,y)$) the sectional curvature of $g$
(resp. $g_{|\dr M}$) on the 2-plane generated by $x,y$. Let $(u_1,
\cdots, u_n)$ be an orthonormal frame on $\dr M$ for which $\II$ is
diagonal, with eigenvalues $k_1, \cdots, k_n$. Then, by the Gauss formula:
$$ \Kb(u_i, u_j) = K(u_i, u_j) + k_i k_j $$
and, taking the trace:
$$ \ricb(u_i, u_i) = \ric(u_i, u_i) - K(u_i, n) + k_i \left( \sum_{j\neq
i} k_j\right)~. $$
Taking the trace once more:
\begin{equation} \label{Sb}
\Sb = S - 2 \ric(n,n) + \sum_i k_i(H-k_i) = S - \frac{2S}{m+1} + H^2 -
\sum_i k_i^2 ~.
\end{equation}

Now
$$ H^2 = \left(\sum_i 1.k_i \right)^2 \leq \left(\sum_i
1^2\right)\left(\sum_i k_i^2\right) = m \sum_i k_i^2$$
with equality if and only if all $k_i$ are equal. The result follows.
$\Box$

\vs

This inequality could be interesting in itself. For instance, if
$S=-m(m+1)$ and $H$ is bounded above by some constant, it implies that
the scalar curvature of $\dr M$ is negative, which has some topological
consequences.

By the way, this computation also leads to the following partial
extension of Theorem \ref{higher}:

\begin{remark}
The second mean curvature $H_2=\sum_{i<j} k_i k_j$ of the boundary of an
Einstein manifold is:
$$ 2 H_2 = \Sb - \frac{m-1}{m+1} S $$
Therefore, it is pointwise constant in an Einstein variation of the
metric which vanishes on the boundary.
\end{remark}

{\bf Proof:} It follows from eq.~(\ref{Sb}). $\Box$

\vs

Now from Proposition \ref{umbilic} and Lemma
\ref{schlafli-einstein} we get:

\begin{cor} \label{cor-umb}
Let $(M, \dr M)$ be a compact manifold with convex (or concave)
boundary. Let $(g_t)_{t_\in
[0,1]}$ be a one-parameter family of Einstein metrics on $M$ with scalar
curvature $S$, such that $\dr M$ is umbilic for $g_0$, and that
the metric induced by $g_t$ on $\dr M$ is constant. Call ${\cal H}_t$ the
integral of the mean curvature of $H$ for $g_t$. Then:
\begin{enumerate}
\item if $S>0$ and $H>0$ (resp. $H<0$) on each connected component of
$\dr M$, then both ${\cal H}_t$ and $V_t$ are minimal (resp. maximal) for
$t_0$, and the variation $H'_t$ of $H_t$ vanishes for $t=0$;
\item if $S<0$ and $H>0$ (resp. $H<0$) on each connected component of
$\dr M$, then ${\cal H}_t$ is minimal and $V_t$ is maximal
(resp. ${\cal H}_t$ is maximal and $V_t$ is minimal) for
$t_0$, and the variation $H'_t$ of $H_t$ vanishes for $t=0$;
\item if $S=0$, then $\dr M$ is umbilic for all $t\in [0,1]$,
i.e. its second fundamental form does not change.
\end{enumerate}
\end{cor}

{\bf Proof:} By Proposition \ref{umbilic}, $H^2$ is pointwise
minimal over $\dr M$ when $\dr M$ is umbilic. If, for instance,
$H>0$ on each connected component of $\dr M$, this shows that
${\cal H}_t$ is also minimal for $t=0$. By Corollary \ref{cor-e},
$V$ is also minimal for $t=0$ when $S>0$, and maximal when $S<0$.
This proves assertions 1 and 2.

For assertion 3, the integral of $H$ over $\dr M$ is constant by
Corollary \ref{cor-e}, while $H^2$ is pointwise minimal over $\dr
M$. Therefore, $H^2$ has to be constant. By the equality case in
Proposition \ref{umbilic}, $\dr M$ has to remain umbilic in that
case. $\Box$

The Corollaries \ref{extension} and \ref{cor-umb} lead to a
description of the non-rigid Ricci-flat manifolds with umbilic
boundary.

\begin{lemma} \label{rigidity}
Suppose $(M, \dr M)$ is a compact $(m+1)$-manifold with boundary, and
$(h_t)_{t\in [0,1]}$ is a non-trivial 1-parameter family of Ricci-flat
metrics on $M$ inducing
the same metric on $\dr M$, and such that $\dr M$ is umbilic for
$h_0$. Then $\dr M$ has at least 2 connected components, and $(h_t)$
corresponds to the displacement of some connected component(s) of $\dr
M$ under the flow of some Killing field(s) of $M$.
\end{lemma}

Note that this is rather restrictive, since a
``generic'' Einstein manifold with boundary should not admit any Killing
field. Some examples are given bellow.

\vs

{\bf Proof:} By Corollary \ref{cor-umb}, $h$ does not change the
induced metric or the second fundamental form of $\dr M$. Call
$\Sigma_1, \cdots, \Sigma_N$ the connected components of $\dr M$.
Then, by Corollary \ref{extension}, each connected component
$\dr_i M$ of $\dr M$ has a neighborhood $\Omega_{i,t}$ which does
not change. Therefore, the deformation $(h_t)$ corresponds to the
displacement of some $\dr_i M$ under Killing fields on $M$. $\Box$

\vs

We shall give some examples of what can happen; this is easier using the
following elementary:

\begin{prop} \label{examples}
Let $(N, g_0)$ be an Einstein $m$-manifold with scalar curvature $m(m-1) k$,
consider the product $M\times \R$ with the warped metric
$$ g := dt^2 + f(t)^2 g_0 $$
where $f$ is a function defined on some interval $I\subset \R$. Then $g$
is Einstein with
scalar curvature $m(m+1) k'$ if and only if $f''(t)= - k'f(t)$ and
$k=k'f(t)^2 + f'(t)^2$ for all $t$. Then each hypersurface $N\times
\{ t\}$ is umbilic in $M$.
\end{prop}

{\bf Proof:}
First check that the Levi-Civita connection $D$ of $g$ is related to the
Levi-Civita connection $\Db$ of $g_0$ by the following formulas: if $n$
is the unit normal vector to $N\times \{ t\}$, and $x,y$ are vector
fields on $N$ (and their extensions on $M$) then
$$ D_x y = \Db_x y - \frac{f'}{f} g(x,y) n $$
$$ D_n x = D_x n = \frac{f'}{f} x $$
$$ D_n n = 0 $$
This is because those expressions define a torsion-free connection
compatible with $g$.

Now this expressions of $D$ shows that each hypersurface $N\times \{
t_0\}$ in $M$ is umbilic, with second fundamental form
$-(f'/f)g$. Therefore, the Gauss formula shows that the sectional
curvature of $M$ on any 2-plane tangent to $N\times \{
t_0\}$ is:
$$ K_M = \frac{1}{f^2} K_N - \frac{f'^2}{f^2} $$
and on each 2-plane containing the direction normal to $N\times \{
t_0\}$:
$$ K_M = -\frac{f''}{f} $$
which shows that $\ric^M(n,n)= - m f''/f$, and leads to the first
condition.
Taking a trace, we see that for $x$ tangent to $N\times \{ t\}$:
$$ \ric^M(x,x) = \frac{k(m-1) - (m-1) f'^2 + f f''}{f^2} g(x,x) $$
and the second condition follows.
$\Box$

\vs

Since we are interested in Ricci-flat metrics, we have to take $k'=0$,
and we can use this proposition in two ways: either $k=0$ and $f(t)=1$, or
$k=1$ and $f(t)=t$.

The simplest example is the case when $g_0$ is the canonical metric on
the sphere $S^m$:

\begin{ex}
Consider the unit ball $B^{m+1}\subset \R^{m+1}$. Any 1-parameter Einstein
deformation of
the metric in $B^{m+1}$ which doesn't change the induced metric on the
boundary $S^m$  is trivial.
\end{ex}

On the other hand:

\begin{ex}
Consider the ``cylinder'' $\Omega:=T^{m}\times [0,1]$.
There exists a
$1$-parameter family of deformations of the metric on $\Omega$ which
does not change the induced metric on the boundary.
\end{ex}

This deformation is obtained by one of the boundary components along the
axis of the cylinder. This happens because
there is a Killing vector field, which is parallel to the axis.

\section{Codimension one foliations}

We give in this section some simple formulas obtained by applying
Theorem \ref{def-hyp-h} to codimension one foliations of Einstein
manifolds. Let $(\Sigma_t)_{t\in I}$ be a smooth one-parameter
family of hypersurfaces in an Einstein $(m+1)$-manifold with
scalar curvature $m(m+1)K$. Suppose that the $\Sigma_t$ define a
foliation of a domain $\Omega\subset M$. For each $x\in \Omega$,
let $H_2(x)$ be the second mean curvature of $\Sigma_t$ at $x$
(for $t$ such that $x\in \Sigma_t$) and let $S_{\Sigma}(x)$ be the
scalar curvature of $\Sigma_t$ for the induced metric. Then:

\begin{thm} The volume $V(\Omega)$ of $\Omega$ is:
\begin{equation} \label{H2}
m K V(\Omega) = 2\int_{\Omega} H_2 dV + \int_{\dr \Omega} H(\dr \Omega) dA~,
\end{equation}
\begin{equation} \label{HSi}
m K V(\Omega) = \int_{\Omega} H^2_{\Sigma} - \tr(\III_\Sigma)dV +
\int_{\dr \Omega} H(\dr \Omega) dA~,
\end{equation}and also
\begin{equation} \label{S}
m^2 K V(\Omega) = \int_{\Omega} S_{\Sigma} dV + \int_{\dr \Omega} H(\dr
\Omega) dA~.
\end{equation}
\end{thm}

Here $H(\dr \Omega)$ is the mean curvature of $\dr \Omega$. $H_\Sigma$,
$S_\Sigma$ and $\III_\Sigma$ are the mean curvature, sectional curvature
and third fundamental form of $\Sigma_t$.

{\bf Proof: } Suppose for instance that $I=[0,1]$, and  note: $$
V_t = \mbox{Vol}\left( \cup_{s=0}^t \Sigma_s \right)~. $$ Denote
again the variations of $V_t, I_t$ and $H_t$ by a prime. Choose a
parametrization $\phi_t:\Sigma\rightarrow M$ such that
$\Sigma_t=\phi_t(\Sigma)$, and let $\phi'_t=v+fn$, where $v\in
T\Sigma$ and $n$ is the unit normal vector to $\Sigma_t$. Then, by
eq.~(\ref{schlafli-h}):
\begin{eqnarray}
n K V'_t & = & \int_{\Sigma_t} H_t'+ \frac{1}{2} \langle I'_t,
\II_t\rangle dA \nonumber \\
& = & \frac{d}{dt}\int_{\Sigma_t} H_t dA - \int_{\Sigma_t} H_t dA' +
\int_{\Sigma_t} \frac{1}{2} \langle I'_t, \II_t\rangle dA \nonumber \\
& = & \frac{d}{dt}\int_{\Sigma_t} H_t dA +
\int_{\Sigma_t} \frac{1}{2} \langle I'_t, \II_t\rangle - H_t (-fH_t +
\overline{\mbox{div}}(v)) dA \nonumber \\
& = & \frac{d}{dt}\int_{\Sigma_t} H_t dA +
\int_{\Sigma_t} -f\langle \II_t, \II_t\rangle + \langle \deltab^*v,
\II_t\rangle + fH_t^2 - H_t \overline{\mbox{div}}(v)) dA \nonumber \\
& = & \frac{d}{dt}\int_{\Sigma_t} H_t dA +
\int_{\Sigma_t} 2f H_2
+ dH_t(v) dA \nonumber
\end{eqnarray}
and eq.~(\ref{H2}) follows, because $\deltab\II_t = -dH_t$ (see
(\ref{del})). Eq.~(\ref{HSi}) is a direct consequence because
$H^2=\tr(\III)+2H_2$.

Taking twice the trace of the Gauss equation for $\Sigma_t$ shows
that: $$ S_{\Sigma} = m(m-1) K + 2 H_2~, $$ so that eq.~(\ref{H2})
becomes: $$ m K V(\Omega) = \int_{\Omega} S_{\Sigma} dV - m(m-1) K
V(\Omega) + \int_{\dr \Omega} H dA $$ which proves (\ref{S}).
$\Box$

This leads for instance to the following simple consequence:

\begin{cor} \label{cr-feuil}
If $K>0$, then no open domain of $M$ has a foliation by closed, minimal
hypersurfaces. If $K=0$, any such foliation is by totally geodesic
hypersurfaces.
\end{cor}

{\bf Proof: } Suppose first that $K>0$. Apply eq.~(\ref{HSi}) to
such a foliation. The boundary term vanishes, and the right-hand
side is therefore non-positive, while the left-hand side is
positive, a contradiction.

If $K=0$, the same argument shows that $H_2\equiv 0$, therefore
$\III\equiv 0$ on each hypersurface by eq.~(\ref{HSi}), and each
hypersurface is totally geodesic. $\Box$

\bigskip

This strongly contrasts with the negatively curved case; for instance,
it is conjectured that a hyperbolic 3-manifold which fibers over the
circle admits a foliation by compact minimal surfaces. Equation
(\ref{HSi}) indicates that such a minimal foliation should have a
remarkable property: the Gauss curvature of each leave, integrated
against a weight corresponding to the amplitude of the normal
deformation, should be constant.

Corollary \ref{cr-feuil} is not too difficult to obtain by other
methods; it is interesting to remark, however, that equations (\ref{H2})
and (\ref{HSi}) can also be used to obtain more general results, for
instance to give an integral lower bound on the mean curvature of a
foliation by minimal hypersurfaces in a positively curved Einstein
manifold.

\bigskip

\section{A quick tour of integral geometry}
\label{intgeom}

In this section, we give a summary of some concepts of integral
geometry which permit us to interpret some of our results more
geometrically. Of course, it cannot be hoped that we can give anything
resembling a comprehensive survey. A reader more interested in this
fascinating subject is referred to the treatises of Santal{\'o}
\cite{santalo} and of Burago and Zalgaller \cite{buzal}.

First, we recall some formulas of Crofton type. Consider
$n$-dimensional Euclidean space $E^n$ and consider the Grassmanian
of all the affine $m$-planes in $E^n$ --  $G_m^n$. This has a
measure invariant with respect to the isometry group of $E^n$. Any
two such measures differ by a constant factor. There is a standard
way to normalize, which will be implicit in the identities we
shall state. Anyway, call the ``canonical'' invariant volume form
$dv_m^n$. There is a natural functional defined on (not
necessarily) convex sets $K$ in $E^n$, to wit,
\begin{equation}
P_m(K) = \int_{L_m\cap K \neq \emptyset} dv_m^n.
\end{equation}
Another natural functional is the following: Consider the space of
$m$-dimensional linear subspaces of $E^n$, and consider the average
$m$-dimensional area of the projections of $K$ onto such
subspaces. This is the so-called {\em Quermassintegral} $W_m(K)$.
It is fairly clear that $P_m$ and $W_m$ are related, and indeed, one
of the fundamental formulas of integral geometry
(\cite[eq. (14.1)]{santalo}) is
\begin{equation}
\label{figg}
P_m(K) = \frac{n O_{n-2}O_{n-3}\cdots O_{n-m-1}}{(n-m)O_{m-1}\cdots
  O_1 O_0} W_m(K),
\end{equation}
where $O_i$ is just the surface area of the $i$-dimensional unit
sphere $S^i$ (it should be noted that the fraction in eq. (\ref{figg})
is just the volume of the Grassmanian of the $m$-dimensional subspaces
in $E^n$).

Equation (\ref{figg}) allows us to relate the quantity $P_m(K)$ to
the integral of the $m$-th symmetric function of curvature, as
follows: Consider the volume $V_\epsilon$ of the $\epsilon$
neighborhood of $K$. A simple computation with radii of curvature
shows that if $\kappa = \kappa_1, \dots, \kappa_{n-1}$ is the
vector of principal curvatures of $\partial K$ (we assume that
$\partial K$ is at least $C^2$ smooth), and $\sigma_m(\kappa)$ is
the $m$-th symmetric function of curvature, then
\begin{equation}
\label{nbhdvol}
V_\epsilon(K) = V(K)+ \sum_{m=0}^{n-1} \frac{\epsilon^{m+1}}{m+1}
\int_{\partial K} \sigma_m(K).
\end{equation}

On the other hand, there is another expression for $V_\epsilon(K)$,
  due to Steiner (see \cite[III.13.3]{santalo}):
\begin{equation}
\label{fig2}
V_\epsilon(K) = \sum_{i=0}^n{n\choose i}W_i(K) \epsilon^i.
\end{equation}
Comparing the coefficients of the powers of $\epsilon$ in formulas
  (\ref{nbhdvol}) and (\ref{fig2}), we see that:
\begin{equation}
m \int_{\partial K} \sigma_{m-1}(K)ds = {n \choose m} W_m(K).
\end{equation}

In other words, the integral mean curvatures of $K$ are directly
expressible in terms of the measures of the set of planes
intersecting $K$, and the average projection measures. In
particular, since $\sigma_0$ is equal to $1$, we see that the area
of $\partial K$ is equal to a constant factor times the measure of
the set of lines intersecting $K$, while the total (first) mean
curvature is a constant times the measure of the set of $2$-planes
intersecting $K$. Since we are especially interested in $n=3$, we
will write down the constants explicitly in that case:

$A(\partial K) = 3 W_1(K),$ while $P_1(K) = \frac{6\pi}{4} W_1(K)$, so
\begin{equation}
\label{areaf}
A(\partial K) = \frac{2}{\pi} P_1(K).
\end{equation}

On the other hand,
$\int_{\partial K} (k_1+k_2) = 3 W_2(K),$ while
$P_2(K) = \frac{3}{2} W_2(K),$
so
\begin{equation}
\label{mintf}
\int_{\partial K} (k_1+k_2) = 2 P_2(K).
\end{equation}

So far, we have talked about convex bodies in a Euclidean setting, but
the theory can be extended to other symmetric spaces, in particular,
to $H^n$ and $S^n$. The expressions become a bit more complicated in
general, but in three dimensions, they are simple enough; the
following formulas are in \cite[eq.~(17.62)]{santalo}:

\begin{equation}
\label{fnon}
P_1(K) = \frac{\pi}{2} A(\partial K),\qquad P_2(K) =
\frac{1}{2}\int_{\partial K}(k_1+k_2) + k V(K),
\end{equation}
where $k$ is the sectional curvature of the ambient space (so the
formula reduces to eq.~(\ref{mintf}) when $k=0$). It should be noted
that there is yet another interpretation of the quantity $P_2(K)$, in
terms of the {\em polar map} of a spherical or hyperbolic convex body
(the spherical version is classical, the hyperbolic has been studied
in the first author's thesis \cite{rivhod}): this map associates to
$K$ the set $K^*$ of hyperplanes intersecting $K$. For the sphere
$S^n$, the set $K^*$ can be naturally viewed as a convex body in
$S^n$, for $H^n$ the set of hyperplanes can be naturally viewed as the
de~Sitter space $S_1^{n-1}$

  Alexandrov's inequality (\cite[p.~145]{buzal}) is the following:
  for a convex $K$ in $E^n$, the cross-sectional measures satisfy:
\begin{equation}
\label{alin}
  V_j^i(K)\geq v_n^{i-j}V_i^j(K),\qquad j\geq i,
\end{equation}
with equality if and
only if $K$ is a ball in $R^n$.
Here $V_i$ can be defined by the following formula:
(\cite[p.~140]{buzal}):
\begin{equation}
  V_k = \frac{1}{n {{n-1}\choose{m}}}
\int_{\partial K} \sum_j k_{j_1}   \dots k_{j_{n-m-1}} d F(x),
\end{equation}
where $v_n$ is the volume of the $n$-dimensional unit ball~; $k_i$
($i=1, \dots, n-1$) are the principal curvatures of
$\partial K$ at the point $x$ of $\partial K$, $d F$ is the area
element of $\partial K$, the sum is taken over all possible finite
sequence of indices $j_1, \dots, j_{n-m-1}$. In particular,
$V_{n-2}(K)$ is related to $Q(S)=\int_S H$ by $Q(S) = n(n-1)
V_{n-2}(K).$ Aleksandrov's inequality in $E^3$ can thus be restated
as follows~: Among all convex bodies with a fixed $P_1$, the ball has
the biggest $P_2$. Our Theorem \ref{af} is exactly the extension of
this result to other three-dimensional space forms.

\section{Extending the Aleksandrov inequality}

This section contains applications of the previous results in the simple
setting of three-dimensional space-forms, esp. $H^3$. Thus we consider a
smooth, strictly convex surface $\Sigma$ in a constant curvature space
$M$ which might be $S^3, \R^3, H^3$ or $S^3_1$ (in this case we suppose
that $\Sigma$ is space-like).

To keep notations close to that of the previous section, we define a
functional $P_2$ as:
$$ 2P_2 := \int_\Sigma H da - 2\epsilon K_0 V~, $$
where $K_0$ is the sectional curvature of $M$, and $\epsilon=1$ if $M$ is
Riemannian, $\epsilon=-1$ if $M=S^3_1$.

Note that, for any deformation of $\Sigma$:
\begin{displaymath}
\left(\int_{\Sigma} H da \right)' = \int_{\Sigma} H'da +
\frac{1}{2}\int_{\Sigma} H\langle I', I\rangle da~.
\end{displaymath}
Therefore, as a consequence of equations (\ref{schlafli-h}) and
(\ref{schlafli-ds}), we have for any deformation of $\Sigma$:
\begin{equation} \label{dF-gen}
P_2' = - \frac{1}{4} \int_\Sigma \langle I', \II-HI\rangle da
\end{equation}

Consider a normal deformation of $\Sigma$, i.e. an infinitesimal
deformation by a vector field $f n$, where $n$ is the (exterior) unit
normal to $\Sigma$. Then $I'=2f\II$, so that:
\begin{equation} \label{dF-norm}
2P_2' = \int_\Sigma \langle -f\II, \II-HI\rangle da = \int_\Sigma 2f K_e da
\end{equation}
where $K_e:=\det(\II)$ is the extrinsic curvature of $\Sigma$.
On the other hand, the area of $\Sigma$ varies as:
\begin{equation} \label{dA-norm}
A' = \frac{1}{2}\int_\Sigma \langle 2f\II, I\rangle da = \int_\Sigma
fH da
\end{equation}

As a consequence, we already find an extremely simple proof of a result
with a flavor of classical differential geometry. It can be seen as a
consequence of some more general results (see \cite{ros}, and also
\cite{ecker-huisken}) but we include it here because of its extremely
simple proof.

\begin{thm}
Suppose that $\Sigma$ is a smooth, strictly convex surface in a
3-dimensional space-form, and
that there exists a constant $k\in\R_+$ such that, on $\Sigma$,
$K_e=kH$. Then $\Sigma$ is totally umbilical.
\end{thm}

{\bf Proof:} Suppose that $K_e=kH$. Then, by equations
(\ref{dA-norm}) and (\ref{dF-norm}), $\Sigma$ is a critical point
of $P_2$ among surfaces with the same area. But it is well known
(see \cite{Po}) that all variations $I'$ of $I$ are induced by
deformations of $\Sigma$. Therefore, eq.~(\ref{dF-gen}) shows that
there exists a constant $k'$ such that $\II-HI=k'I$, so that
$\Sigma$ is totally umbilical. $\Box$

\bigskip

We now turn to the extension of the classical Alexandrov
inequality (see Section \ref{intgeom}) for convex surfaces in
Euclidean space to three-dimensional space-forms.

\begin{thm}\label{af}
Let $S$ be a compact convex surface in $H^3$ (resp. $R^3$, $S^3$). Let
$V(S)$ be the volume of the interior of $S$, and call
$2P_2(S):=\int_S H da + 2V(S)$ (resp. $2P_2(S)=\int_S H da$,
$2P_2(S)=\int_S H da-2V(S)$). There exists a (unique modulo global
isometries) umbilical surface
$S_0$ with $\area(S_0)=\area(S)$, and $P_2(S_0)\geq P_2(S)$.
\end{thm}

\begin{remark}\label{af1}
  An elementary consequence is the well-known fact that a convex
  surface in $S^3$ has area at most $2\pi$.
\end{remark}

\bigskip

{\bf Proof. } We will give the proof of Theorem \ref{af} in the
hyperbolic case, the other two situations are very similar. For
$k_0>0$, let ${\cal C}_{A, k_0}$ be the space of smooth, convex
surfaces in $H^3$ with area $A$ and principal curvatures at most
$k_0$, and containing a given fixed point $x_0$. It is again an
elementary consequence of eq.~(\ref{dF-gen}) and of \cite{Po} that
the only critical points of $P_2$ in ${\cal C}_{A, k_0}$ are the
umbilical hypersurfaces (there is a unique such surface in ${\cal
C}_{A, k_0}$ modulo the global isometries of $H^3$).

It is therefore sufficient to prove that $P_2$ has a maximum in
the interior of ${\cal C}_{A, k_0}$. This will be a consequence of the
following points:
\begin{enumerate}
\item if $S\in {\cal C}_{A, k_0}$ contains a point $s$ where a principal
curvature vanishes, then there exists an infinitesimal deformation of
$S$ increasing
the minimum of the principal curvatures and increasing $P_2$, while
leaving the area constant;
\item if $k_0>1$ and if $S\in {\cal C}_{A, k_0}$ contains a point $s$
where a principal
curvature is equal to $k_0$, then there exists a deformation of $S$
decreasing
the maximum of the principal curvatures and increasing $P_2$;
\item for each $M>0$, there exists $L>0$ such that if $S\in {\cal C}_{A,
k_0}$ has (extrinsic) diameter above $L$, then $P_2(S)\leq -M$.
\end{enumerate}

Theorem \ref{af} follows, because a maximizing sequence for $P_2$ can
neither ``degenerate'' (because of point (3)), nor converge to a surface
with a vanishing principal curvature (point (1)) or a principal
curvature equal to $k_0$ (point (2)).

To prove point (1), note that the equation (\ref{dII-norm})
simplifies here to:
\begin{equation} \label{dIIb}
\II' = H_f + fI - f \III
\end{equation}
Therefore, to insure that a deformation $f n$ increases the
minimum of the principal curvatures at a point where this minimum
vanishes, it is enough to have: $f\geq \| H_f\|$ over $S$. But, if
a principal curvature of $S$ vanishes, then $S$ is not umbilical,
so that $H$ and $K_e$ are not proportional on $S$. It is then easy
to check using equations (\ref{dF-norm}) and (\ref{dA-norm}) that
there exists a normal deformation of $S$ with the right
properties.

Point (2) can be proved in the same way, the condition on $f$ is now
that $(k_0-1)f\geq \| H_f\|$.

Finally, for point (3), let $S$ be a convex surface in $H^3$, and, for
$\epsilon\geq 0$, call
$E_\epsilon$ the set of points at distance at most $\epsilon$ from the
interior of $S$, and $S_\epsilon:=\dr E_\epsilon$. Let ${\cal
H}_\epsilon$ be the integral mean curvature of $S_\epsilon$, and let
${\cal A}_\epsilon$ be its area and ${\cal V}_\epsilon$ the volume of
its interior $E_\epsilon$. Equation
(\ref{schlafli-h}) shows that:
$$ -2 \frac{d {\cal V}_\epsilon}{d\epsilon} = \frac{d {\cal
H}_\epsilon}{d\epsilon} + \int_{S_\epsilon} \langle \II, \II-HI\rangle
da~, $$
so that:
$$ 2 \frac{d {\cal V}_\epsilon}{d\epsilon} + \frac{d {\cal
H}_\epsilon}{d\epsilon} = \int_{S_\epsilon} K_e da = \int_{S_\epsilon}
(K+1) da = 4\pi + {\cal A}_\epsilon~, $$
where we have used the Gauss-Bonnet theorem. Since $d{\cal
V}_\epsilon/d\epsilon = {\cal A}_\epsilon$, we have:
$$ \frac{d {\cal V}_\epsilon}{d\epsilon} + \frac{d {\cal
H}_\epsilon}{d\epsilon} = 4\pi~, $$
and therefore:
$$ \frac{d^3 {\cal V}_\epsilon}{d\epsilon^3} + \frac{d {\cal
V}_\epsilon}{d\epsilon} = 4\pi~. $$
Integrating this EDO leads to a classical formula for ${\cal
V}_\epsilon$:
\begin{equation} \label{eqvol}
{\cal V}_\epsilon = {\cal A}_0 \sinh(\epsilon) + 4\pi
(\epsilon-\sinh(\epsilon)) +
{\cal H}_0 (\cosh(\epsilon)-1) + {\cal V}_0
\end{equation}
Now suppose that $S$ has extrinsic diameter at least $L$, then
$E_0$ contains a segment $\gamma$ of length at least $L$. Applying
equation (\ref{eqvol}) to a sequence of convex surfaces in $E_0$
which converges to $\gamma$, we find a lower bound for ${\cal
V}_\epsilon$: $$ {\cal V}_\epsilon \geq 4\pi
(\epsilon-\sinh(\epsilon)) + 2\pi L (\cosh(\epsilon)-1) $$ so
that, for any $\epsilon\geq 0$: $$ {\cal A}_0 \sinh(\epsilon)
+{\cal H}_0 (\cosh(\epsilon)-1) + {\cal V}_0 \geq 2\pi L
(\cosh(\epsilon)-1) $$ Taking the limit as $\epsilon\rightarrow
\infty$ shows that: $$ {\cal A}_0 + {\cal H}_0 \geq 2\pi L $$ and
point (3) follows.

This finishes the proof of Theorem \ref{af}. $\Box$

\bigskip

Note that the ``local'' part of this argument extends partly to higher
dimensions, again for Einstein manifolds with convex boundaries. Namely,
if $(M, \dr M)$ is such a manifold, we say that it is ``rigid'' if it
admits no infinitesimal Einstein deformation which does not change the
induced metric on the boundary. It is proved in \cite{ecb} that, in that
case, all infinitesimal deformations of the metric on $\dr M$ are
induced (uniquely) by Einstein deformations of $M$. Therefore, it is
still true in that case that, if $(M, \dr M)$ is a critical point of $P_2$
among Einstein metrics with the same ``area'', then it is umbilical. We
do not know whether there exists any non-rigid Einstein metric with
negative curvature and strictly convex boundary.

\section*{Acknowledgements} The authors are happy to acknowledge their
debt to Fred Almgren, without whom this paper would not have been
possible. I.~Rivin would like to thank the Institut Henri Poincar{\'e}
and the Institut des Hautes {\'E}tudes Scientifiques for their
hospitality at crucial moments.

\bibliographystyle{alpha}

\end{document}